\def\mystyle{}

\if\mystyle l
\documentclass[landscape]{amsart}
\else
\documentclass{amsart}
\fi

\usepackage{begnac}

\newtheorem*{exc*}{Exercise}

\begin{document}

\title{On the expressive power of quantifiers in continuous logic}

\author{Itaï \textsc{Ben Yaacov}}

\address{Itaï \textsc{Ben Yaacov},
  Univ Lyon,
  Université Claude Bernard Lyon 1,
  Institut Camille Jordan, CNRS UMR 5208,
  43 boulevard du 11 novembre 1918,
  69622 Villeurbanne Cedex,
  France}

\urladdr{\url{http://math.univ-lyon1.fr/~begnac/}}

\thanks{Author supported by ANR project AGRUME (ANR-17-CE40-0026).}
\thanks{The author wishes to thank Todor \textsc{Tsaknov} for having pushed him into writing this note, as well as for some mushy chocolate}

% \date{\today}

\begin{abstract}
  In this short note we compare the expressive power of \emph{real-valued continuous logic} (or just \emph{continuous logic}, in recent literature) with that of \emph{compact-valued continuous logic}, proposed by Chang and Keisler.
  We conclude that the two logics have the same expressive power, and moreover, that this remains true if we replace the plethora of potential quantifiers of compact-valued logic with a single ``primordial'' one.
\end{abstract}

\maketitle

% \tableofcontents

This note summarises remarks the author has made on various occasions, and to various people, regarding the similarities and differences between two variants of ``continuous logic'' that exist in the literature, and how to reconcile them: the one proposed by Chang and Keisler \cite{Chang-Keisler:ContinuousModelTheory}, which we shall call \emph{compact-valued} logic, and that proposed in \cite{BenYaacov-Usvyatsov:CFO}, which we call here \emph{real-valued}.
Despite striking similarities, the latter was conceived in ignorance of the former (at least as far as the author of the present note is concerned), essentially through reverse-engineering the idea of compact type spaces with open continuous variable restriction maps between them.
Let us also refer to the \textit{True}/\textit{False} logic one learns in kindergarten as \emph{two-valued}.

For most of this note we are going to consider ``bare bones'' versions of either logic, in a purely relational language, without a distinguished symbol for equality or distance.
The latter means that a structure may contain two points that are formally distinct and yet indistinguishable by the logic, even with parameters.
For our purposes this will be of no importance whatsoever (see a further remark to this effect at the end).
Relation symbols applied to variable symbols (in the absence of function symbols) make atomic formulas, and these are combined using continuous connectives and appropriate first order quantifiers to generate all the formulas of the logic.

In compact-valued logic, the truth values of a formula $\varphi$ belong to a fixed compact value space (this note being written in France, compact spaces are Hausdorff).
In particular, the logic associates to each relation symbol $P$ a compact value space $X$, and then the atomic formula $P(x)$ is $X$-valued.
If $(\varphi_i : i < \alpha)$ is a sequence of, respectively, $X_i$-valued formulas, and $\theta\colon \prod X_i \rightarrow Y$ is continuous, then we view $\theta$ as a \emph{connective} that combines the given formulas into a $Y$-valued formula $\theta \circ (\varphi_i : i < \alpha)$.

In real-valued logic, formulas take values in compact intervals of the reals, and connectives are continuous functions $\theta\colon \bR^\alpha \rightarrow \bR$ (necessarily bounded on each product of compact intervals).
In some texts, one restricts truth values to $[0,1]$, but the result is the same.

\begin{rmk*}
  With these definitions, real-valued formulas are closed under limits of uniformly Cauchy sequences.
  Indeed, for any such a sequence $(\varphi_n)$ there exists a continuous $\theta\colon \bR^\bN \rightarrow \bR$ such that $\theta \circ (\varphi_n)$ is the limit.
  Moreover, $\theta$ can be chosen in a manner that depends only on the convergence rate of $(\varphi_n)$, so, possibly passing to a sub-sequence, a single such $\theta$ suffices (see the \emph{forced limit} construction in \cite{BenYaacov-Usvyatsov:CFO}).
  By the Stone-Weierstraß Theorem, closing under this special $\theta$ plus, say, constants, addition and multiplication, is equivalent to closing under all possible continuous connectives.
\end{rmk*}

Clearly, atomic real-valued formulas are also atomic compact-valued ones, and real-valued continuous connectives, once restricted to an appropriate product of compact intervals, are also compact-valued ones.

In the opposite direction, consider a compact-valued language $L$.
For each $X$-valued relation symbol $P$ of $L$, choose an embedding $X \subseteq [0,1]^\alpha$, and let $(P_i : i < \alpha)$ be new $[0,1]$-valued relation symbols of the same arity.
Let $L^\bR$ be the real-valued language consisting of all such symbols.
An $L$-structure $M$ will be tacitly identified with an $L^\bR$-structure $M^\bR$ defined so $P(a) = \bigl( P_i(a) : i < \alpha \bigr) \in X \subseteq [0,1]^\alpha$.

Say that an $X$-valued $L$-formula $\varphi(x)$ is \emph{coded} in $L^\bR$ if for every continuous function $\theta\colon X \rightarrow \bR$, there exists an $L^\bR$-formula $\varphi_\theta(x)$ that agrees with $\theta \circ \varphi$ on all $L$-structures.
An equivalent condition is that the same hold for a family $\Xi$ of functions $\xi\colon X \rightarrow \bR$ that separates points in $X$.
Indeed, under this assumption, the induced map $X \rightarrow \bR^\Xi$ is continuous and injective, and therefore a topological embedding.
By Tietze's Extension theorem, any continuous $\theta\colon X \rightarrow \bR$ extends continuously to $\bR^\Xi$, and then $\varphi_\theta = \theta \circ (\varphi_\xi : \xi \in \Xi)$ will do.

By construction, every atomic $L$-formula is coded in $L^\bR$.
By the same Tietze argument as above, any continuous combination of coded formulas is again coded.
It follows that every quantifier-free $X$-valued $L$-formula is coded in $L^\bR$.
In other words, at the quantifier-free level, the expressive power of real-valued and compact-valued logic is the same.

For real-valued logic we proposed $\sup$ and $\inf$ as quantifiers, which requires some justification.
A first argument is that, \textit{a priori}, these may be presented as the real-valued analogues of the classical quantifiers $\forall$ and $\exists$.
A second argument, \emph{a posteriori}, is that these quantifiers ``deliver the goods'', in the sense that real-valued logic does behave as a generalisation of two-valued logic.
In particular, when stated appropriately, Łoś's Theorem and the Compactness Theorem hold in real-valued logic, as do many other results.
The discussion which follows will provide, in a sense, a third argument, which may be qualified as \textit{ab initio}: the choice of $\sup$ and $\inf$ as quantifiers is correct because they have the same expressive power as the \emph{primordial quantifier}, from which any (admissible) quantifier may be constructed.

In a general compact-valued logic there are no obvious analogues for $\forall$ and $\exists$.
Chang and Keisler proposed that quantifiers should be specified at the same time as other symbols of the language.
A compact-valued quantifier $Q$, in this sense, is a function that takes a collection of truth values in some value space $X$ (e.g., all values of $\varphi(a,b)$ in a structure $M$, as $a$ varies over $M$) and returns a truth value in some value space $Y$ (the value of $(Qx) \varphi(x,b)$ in $M$).
Notice that $Qx$ can only be applied to a formula $\varphi(x,y)$ whose value space is $X$, and then the value space of $(Qx) \varphi(x,y)$ is $Y$.
In order to be admissible as a quantifier, $Q$ is also required to be continuous in a sense that we shall make precise below.
With such continuous quantifiers, compact-valued logic is well-behaved (in particular, Łoś and Compactness hold).

\begin{dfn*}
  Let $X$ be a Hausdorff topological space, and let $\cK(X)$ denote the space of all non-empty compact subsets of $X$.
  For an open set $U \subseteq X$, let
  \begin{gather*}
    O^1_U = \bigl\{ K \in \cK(X) : K \subseteq U \bigr\},
    \qquad
    O^2_U = \bigl\{ K \in \cK(X) : K \cap U \neq \emptyset \bigr\}.
  \end{gather*}
  We equip $\cK(X)$ with the \emph{Vietoris Topology}, namely the topology generated by all sets $O^1_U$ and $O^2_U$, as $U$ varies over all open subsets of $X$.
\end{dfn*}

A basic open set of the Vietoris topology is of the form $O^1_U \cap \bigcap_{i<n} O^2_{V_i}$.
It is Hausdorff, and if $X$ is compact, then so is $\cK(X)$.
See also Kechris \cite[Chapter~I.4.F]{Kechris:Classical}, where $\cK(X)$ also contains the isolated point $\emptyset$.

We can now state the continuity requirement for a compact-valued quantifier: it should be of the form $Q(A) = \theta_Q(\overline{A})$, where $\theta_Q\colon \cK(X) \rightarrow Y$ is continuous (and $\emptyset \neq A \subseteq X$, so $\overline{A} \in \cK(X)$).

\begin{exm*}
  If $I \subseteq \bR$ is a compact interval, then the upper bound function $\sup \colon \cK(I) \rightarrow I$ is continuous in the Vietoris topology, and similarly for $\inf$.
  Therefore, the real-valued quantifiers are also compact-valued quantifiers in the sense of Chang and Keisler.
\end{exm*}

Let us simplify things and define a variant of compact-valued logic in which the aforementioned functions $\theta_Q\colon \cK(X) \rightarrow Y$ are exactly all possible identity maps.
Stated somewhat differently, we allow a unique quantifier $Q$, which may be applied to any formula, with any value space.
If $\varphi(x,y)$ is an $X$-valued formula, then $(Qx) \varphi(x,y)$ is $\cK(X)$-valued.
If $M$ is a structure and $b \in M$, then, with truth values calculated in $M$:
\begin{gather*}
  (Qx) \varphi(x,b) = \overline{\bigl\{ \varphi(a,b) : a \in M \bigr\}} \in \cK(X).
\end{gather*}
Let us call this $Q$ the \emph{primordial quantifier}.
For any other continuous quantifier $Q'$, in the sense of Chang and Kiesler, we may express $(Q'x)\varphi(x,y)$ as a composition $\theta_{Q'} \circ (Qx) \varphi(x,y)$, now viewing $\theta_{Q'}$ as a mere continuous connective.

Let us compare he expressive power of quantifiers in real-valued and in compact-valued logic.
We have already observed that the real-valued quantifiers $\sup$ and $\inf$ can be expressed as a composition of continuous functions with the primordial quantifier.
For the opposite direction, let us observe first that $\cK$ is a functor.
Indeed, if $\theta\colon X \rightarrow Y$ is continuous, and $K \subseteq X$ is compact, then $\theta(K) \subseteq Y$ is compact as well, defining a map $\cK(\theta) \colon \cK(X) \rightarrow \cK(Y)$.
In addition, assuming that $U \subseteq Y$ is open and $V \subseteq X$ is its inverse image by $\theta$, then $O^1_V$ and $O^2_V$ are the inverse images by $\cK(\theta)$ of $O^1_U$ and $O^2_U$, respectively, so $\cK(\theta)$ is continuous.
In particular, if $\theta\colon X \rightarrow [0,1]$ is continuous, then so is the function $\sup \circ \cK(\theta) \colon \cK(X) \rightarrow [0,1]$, which we shall denote by $(\sup \theta)$.
Assume that $K, F \in \cK(X)$ are distinct, say $x \in F \setminus K$.
By Urysohn's Lemma there exists a continuous $\theta\colon X \rightarrow [0,1]$ that vanishes on $K$ and equals one at $x$.
Then $(\sup \theta)$ vanishes at $K$ and equals one at $F$.
It follows that the functions $(\sup \theta)$ separate points in $\cK(X)$.

Assume now that $\varphi$ is an $X$-valued $L$-formula, in compact-valued logic, and that it is coded in $L^\bR$.
In other word, for every continuous $\theta\colon X \rightarrow [0,1]$, the $L$-formula $\theta \circ \varphi$ agrees with a real-valued $L^\bR$-formula $\varphi_\theta$.
Then $(\sup \theta) \circ (Qx)\varphi$ agrees with $\sup_x \varphi_\theta$, and since the functions $(\sup \theta)$ separate points, $(Qx)\varphi$ is also coded in $L^\bR$.

We conclude that real-valued continuous logic, equipped with the two quantifiers $\sup$ and $\inf$ (or just with one, from which the other can be recovered), has exactly the same expressive power as the compact-valued logic of Chang and Keisler, even when restricting the latter to the primordial quantifier alone.

In order to make this a little more precise, define a \emph{compact-valued theory} $T$ to be (the set of consequences of) a collection of \emph{closed conditions} $\sigma \in K$, where $\sigma$ is, say, an $X_\sigma$-valued sentence, and $K \subseteq X_\sigma$ is closed.
We could define a real-valued theory in the same manner, but we usually prefer to restrict to sets of conditions of the form $\sigma = 0$, which may then be identified with just a collection of sentences.
With these definitions, the class of $L^\bR$-structures that code $L$-structures is elementary, defined by a (universal) \emph{base theory} $T_0$.
General classes of models of $L$-theories are in bijection with classes of models of $L^\bR$-theories extending $T_0$.

\begin{rmk*}
  Two-valued logic usually assumes the presence of a distinguished relation symbol that interprets equality.
  When constructing real-valued logic, one argues that the ``natural'' real-valued analogue of equality is a symbol that interprets a (bounded, complete) distance function, with respect to which all other symbols are uniformly continuous.
  In addition to being natural, this is extremely useful for applications, making real-valued logic suitable for the model-theoretic treatment of metric structures.
  This potential was unfortunately missed by Chang and Keisler, who followed the tradition of two-valued logic, assuming an equality symbol.

  In order to designate a real-valued binary relation symbol $d$ as a distance symbol, we need to add to our base theory the universal axioms asserting that $d$ interprets a pseudo-distance with respect to which every other symbol respects a fixed uniform continuity modulus (restricting the value space of $d$ to $\{0,1\}$ would make it an equality symbol).
  If $M$ is a model of these axioms then we may quotient out the zero-distance relation, and complete it (as a structure), without affecting the interpretation of the formulas.
\end{rmk*}

\begin{rmk*}
  Once we have a distance (or equality) symbol, we may introduce function symbols (without one, function symbols may misbehave).
  Adding a (uniformly continuous) function symbol $f(x)$ is the same as adding a relation symbol $P(x,y)$ representing $d\bigl( f(x), y \bigr)$, together with $\forall\exists$ axioms that assert that $P$ is indeed of that form.
\end{rmk*}

\bibliographystyle{begnac}
\bibliography{begnac}

\end{document}